\documentclass[12pt]{amsart}
\usepackage{amsfonts, amsmath, amstext, amssymb, amsthm}
\usepackage{calrsfs}
\usepackage{latexsym}
\usepackage{fullpage}

\newcommand{\calH}{{\mathcal H}}

\newcommand{\calM}{{\mathcal M}}

\newcommand{\calO}{{\mathcal O}}
\newcommand{\calP}{{\mathcal P}}

\newcommand{\calR}{{\mathcal R}}

\newcommand{\calT}{{\mathcal T}}
\newcommand{\calU}{{\mathcal U}}
\newcommand{\calV}{{\mathcal V}}
\newcommand{\calW}{{\mathcal W}}

\newcommand{\Z}{{\mathbb Z}}

\newcommand{\ord}{\mbox{ord}}
\newcommand{\pp}{{\mathfrak p}}

\newcommand{\calS}{{\mathcal S}}

\newcommand{\qq}{{\mathfrak q}}

\newtheorem{theorem}{Theorem}[section]
\newtheorem*{theorem*}{Theorem}
\newtheorem{prop}[theorem]{Proposition}
\newtheorem{corollary}[theorem]{Corollary}

\newtheorem{question}[theorem]{Question}
\newtheorem{lemma}[theorem]{Lemma}

\theoremstyle{definition}
\newtheorem{definition}[theorem]{Definition}

\theoremstyle{remark}
\newtheorem{remark}[theorem]{\bf Remark}
\newtheorem{notation}[theorem]{Notation}

\usepackage{color}

\newcommand{\set}[2]{\ensuremath{ \{ #1 : #2 \} }}

\newcommand{\N}{\mathbb{N}}
\newcommand{\Q}{\mathbb{Q}}

\newcommand{\xvec}{\vec{x}}

\newcommand{\HTP}{\operatorname{HTP}}

\newcommand{\la}{\langle}
\newcommand{\ra}{\rangle}

\newcommand{\comment}[1]{}

\newcommand{\rg}[1]{\text{range}(#1)}

\def\bfz{\boldsymbol{0}}
\def\s01{\ensuremath{\Sigma^0_1}}
\def\d02{\ensuremath{\Delta^0_2}}
\def\phi{\varphi}

\begin{document}

\title{As Easy as $\Q$:  Hilbert's Tenth Problem\\for Subrings of the Rationals and Number Fields}
\author{Kirsten Eisentr\"ager} \thanks{This project was initiated at a workshop held at
the American Institute of Mathematics, after an extended conversation with Tom Scanlon.\\
\indent
The first author was partially supported by National
  Science Foundation grant DMS-1056703}
\curraddr{Department of Mathematics, The Pennsylvania State
  University, University, Park, PA 16802, USA}
\email{eisentra@math.psu.edu}
\urladdr{http://www.personal.psu.edu/kxe8/}
\author{Russell Miller}\thanks{The second author was partially supported
by N.S.F.\ grants DMS-1001306 and DMS-1362206
, and by PSC-CUNY Research Awards 67839-00 45 and
66582-00 44 from the City University of New York .} 
\curraddr{Dept.\ of Mathematics, Queens College, Queens, NY 11367 \& Ph.D.\ Programs
in Mathematics and Computer Science, CUNY Graduate Center, New York, NY 10016, USA}
\email{Russell.Miller@qc.cuny.edu}
\urladdr{http://qcpages.qc.cuny.edu/~rmiller}
\author{Jennifer Park}\thanks{The third author was partially supported by N.S.F.\ grant
DMS-1069236, and by an NSERC PDF grant.}
\curraddr{Department of Mathematics, McGill University, Montr\'eal, Qu\'ebec, Canada}
\email{jennifer.park2@mcgill.ca}
\urladdr{http://www.math.mcgill.ca/jpark/}
\author{Alexandra Shlapentokh} \thanks{The fourth author was partially
supported by N.S.F.\ grant DMS-1161456}
\curraddr{Department of Mathematics, East Carolina University, Greenville, NC 27858, USA}
\email{shlapentokha@ecu.edu}
\urladdr{}
\maketitle
\begin{abstract}
  Hilbert's Tenth Problem over the rationals is one of the biggest open problems in the area of  undecidability in number theory. In this
  paper we construct new, computably presentable subrings $R \subseteq \Q$
  having the property that Hilbert's Tenth Problem for $R$, denoted $\HTP(R)$,
  is Turing equivalent to $\HTP(\Q)$.

  We are able to put several additional constraints on the rings $R$
  that we construct.  Given any computable nonnegative real number $r\leq 1$ we
  construct such rings $R=\Z[\calS^{-1}]$ with $\calS$ a set of
  primes of lower density $r$. We also construct examples of rings
  $R$ for which deciding membership in $R$ is Turing equivalent to
  deciding $\HTP(R)$ and also equivalent to deciding $\HTP(\Q)$.
  Alternatively, we can make $\HTP(R)$ have arbitrary computably enumerable degree
  above $\HTP(\Q)$.  Finally, we show that the same can be done
  for subrings of number fields and their prime ideals.
\end{abstract}

\section{Introduction}
\label{sec:intro}

Hilbert's Tenth Problem asks to find an algorithm that takes as input 
any polynomial $f \in \Z[x_1, \ldots, x_n]$ and decides
whether $f = 0$ has a solution in $\Z^n$.
In 1969, Matiyasevich~\cite{Mat70}, using
work by Davis, Putnam and Robinson~\cite{DPR61}, proved that no
such algorithm exists; thus, we say that Hilbert's Tenth Problem is
\textit{undecidable}.  We also ask the same question for polynomial equations with coefficients
and solutions in other \textit{computably presentable} rings $R$.  We call this {\it Hilbert's Tenth Problem  over $R$}.

As yet, the answer for $R =\Q$ is unknown.
However, if $\Z$ admits a diophantine definition over $R$, or more
generally, if there is a diophantine model of the ring $\Z$ over $R$,
then a negative answer for Hilbert's Tenth Problem over $R$ can be
deduced from that of $\Z$. This is the reason for the following
definitions below:

\begin{definition}
\label{D:diophantine definition}
Let $R$ be a commutative ring. 
Suppose $A \subseteq R^k$ for some $k \in \N$.
Then $A$ is {\em diophantine over $R$}
if there exists a polynomial $f$ in $k+n$ variables with coefficients
in $R$ such that
\[
    A=\{
{\mathbf t} \in R^k: \,    \exists\, x_1,\ldots,x_n \in R, \; f({\mathbf t},x_1,...,x_n) = 0\}.
\]
\end{definition}


\begin{definition}
\label{D:diophantine model}
A {\em diophantine model of $\Z$ over a ring $R$} is a
subset $A \subseteq R^k$ (for some $k >0$) that is diophantine over $R$,
together with a bijection $\Z \to A$,
under which the graphs of addition and multiplication
(which are subsets of $\Z^3$) correspond to subsets
of $A^3 \subseteq R^{3k}$ that are diophantine over $R$.
\end{definition}

In 1992, Mazur~\cite{M1} conjectured that if $X$ is a variety over
$\Q$, then the topological closure of $X(\Q)$ in $X(\mathbb{R})$ has
finitely many components. This implies (see \cite{CZ}) that $\Z$ is
not diophantine over $\Q$, and that there is no diophantine model of
$\Z$ over $\Q$.  The possible lack of diophantine models and definitions
of $\Z$ over $\Q$ left no obvious way to approach the problem of
showing that $\HTP(\Q)$ is undecidable and motivated the search for
potentially easier problems that might shed some
light on the diophantine problem of $\Q$.  One study of partial
problems led to a reduction of the number of universal quantifiers in
the first-order definition of $\Z$ over $\Q$ and culminated in the
result of Koenigsmann.  Another path led to an exploration of the
rings between $\Z$ and $\Q$ to see if any of them admit a diophantine
definition or a diophantine model of the integers. Subrings of
$\mathbb{Q}$ are in bijection with subsets $\calS$ of the set
$\mathcal P$ of the prime numbers; one associates to $\calS$ the ring
$\Z[\mathcal S^{-1}]$.

If $\calS$ is finite, one can obtain a diophantine
definition of $\Z$ over rings $\Z[\calS^{-1}]$ from~\cite{Rob1},
and hence $\HTP(\Z[\calS^{-1}])$ is undecidable.
From the same work of Robinson~\cite{Rob1}, it also follows that if $\calP-\calS$
is finite, then $\Z[\calS^{-1}]$ is diophantine over $\Q$.  Therefore,
for such $\calS$,
$\HTP(\Z[\calS^{-1}])$ is decidable if and only if $\HTP(\Q)$ is
decidable.    

However, when $\calS$ is both infinite and co-infinite (i.e.\
$\calP - \calS$ is infinite), the results concerning $\Z[\calS^{-1}]$
are more interesting: the first result was obtained by Poonen:
\begin{theorem*}[\cite{Po2}]
There exist computable
sets $\calS$ of natural density zero and of natural density one such
that $\HTP(\Z[\calS^{-1}])$ is undecidable.
\end{theorem*}
 This remarkable paper was
followed by generalizations in \cite{ES}, \cite{Perlega}, and
\cite{EES}. However, no attempt has been made so far in trying to compare $\HTP(\Z[\calS^{-1}])$ to $\HTP(\Q)$. As it is yet unknown whether $\HTP(\Q)$ is decidable or not, one could try to compare the difficulties of $\HTP(\Z[\calS^{-1}])$
for different $\calS$, using the notion of Turing reducibility:

\begin{definition}
Given a set $B \subseteq R$, an {\em oracle for $B$} takes as input an element of $R$, and outputs YES or NO, depending on whether the element belongs to $B$. For $A,B \subseteq R$, $A$ is {\em Turing reducible to $B$} (written $A \leq_T B$) if there is an
algorithm that determines membership in $A$ using an
oracle for $B$. We say that {\em $A$ is Turing equivalent to $B$} ($A \equiv_T
B$) if $A \leq_T B$ and $B \leq_T A$. The set of equivalence classes under $\equiv_T$ are called {\em Turing degrees}.  (The notion of an oracle is also discussed in Section \ref{S: definitions}.)
\end{definition}

The first result of this paper concerns the subrings $R$ of $\Q$ for which
$\HTP(R) \equiv_T \HTP(\Q)$:

\begin{theorem*}[Theorem \ref{thm:low}]
  For every computable real number $r$ between 0 and 1 there exists a
  computably enumerable set $\calS$ of primes of lower density $r$ such that $\HTP(\Z[\calS^{-1}]) \equiv_T  \HTP(\Q)$.
 \end{theorem*}
 
We have imposed the condition that $\calS$ is \textit{computably enumerable} because of the following reason: since Hilbert's Tenth Problem is primarily concerned with algorithms for rings, we want to consider rings in which addition and multiplication are computable, called computably presentable rings (see \S \ref{S: definitions} for a precise definition). 
 
For $\calS \subset \calP$, if the ring $R = \Z[\calS^{-1}]$ is computably presentable, then there is an algorithm such that, when left running forever, prints out exactly $\calS$. In this case, $\calS$ is said to be computably enumerable (see \S \ref{S: definitions} for a precise definition). Conversely, if $\calS$ is computably enumerable, then $R$ is computably presentable.

Matiyasevich showed that every computably enumerable subset of the integers is
diophantine. This implies that $\HTP(\Z)$ is
Turing equivalent to the halting set.  Since there are infinitely many
different Turing degrees of undecidable computably enumerable sets (constructed by \cite{Friedberg} and \cite{Muchnik},
who independently invented the priority method), this is much
stronger than showing that Hilbert's Tenth Problem over $\Z$ is
undecidable. 

To our knowledge, all attempts so far to prove undecidability for $\HTP(\Q)$ 
try to prove that $\HTP(\Q) \equiv_T \HTP(\Z)$. However, it is still possible that $\HTP(\Q) <_T \HTP(\Z)$ while still being undecidable.   One could argue that the results of this paper do not point in this direction.\\

In \S \ref{sec:HTPQoracle}, we construct a set $\calS$ that is not necessarily
computably enumerable (so $\Z[\calS^{-1}]$ is not necessarily computably presentable),
but for which $\HTP(\Z[\calS^{-1}])\equiv_T\HTP(\Q)$.
Then in \S \ref{sec:computable}, we modify these constructions
to make $\calS$ computably enumerable, thus making $\Z[\calS^{-1}]$ computably presentable; this process is more technical, using priority
constructions from computability theory. We hope that the initial
proofs will afford the reader some insight before encountering the details.
We also get an infinite sequence of nested rings $R_i$ with
$\HTP(R_i) \equiv_T \HTP(\Q)$ (see Corollary \ref{cor:nested}).
 
The next theorem produces examples of a different flavour.  Let $\calS$
be a computably enumerable set of primes and let $R=
\Z[\calS^{-1}]$. We trivially have $\HTP(R)$ $\geq_T R$, but it is possible to get $\HTP(R) \equiv_T R$, while also
choosing $\calS$ such that $\HTP(R) \equiv_T \HTP(\Q)$:

\begin{theorem*}[Special case of Theorem \ref{thm:B}]
  There exists a computably presentable ring
  $R=\Z[\calS^{-1}]$, with $\calS$ a computably enumerable subset of
  the prime numbers of lower density 0, such that $R
  \equiv_T \HTP(R)\equiv_T \HTP(\Q)$.
\end{theorem*}
The general version of Theorem \ref{thm:B} shows that many of the preceding
arguments can in fact produce computably presentable subrings $R\subseteq\Q$
for which $\HTP(R)$ is Turing-equivalent to an arbitrary c.e.\ set $S$, provided only
that $S\geq_T\HTP(\Q)$.

The final two results are the complementary ring versions as in \cite{ES} and \cite{EES}.
\begin{theorem*}[Theorem \ref{thm:complements}]
  For any positive integer $m$, the set of all rational primes $\calP$
  can be represented as a union of pairwise disjoint sets $\calS_1
  \ldots, \calS_{m}$, each of upper density 1 and such that for all
  $i$ we have that $\HTP(\Z[\calS_i^{-1}] )\leq_T \HTP(\Q)$ and
  $\calS_i \leq_T \HTP(\Q)$.
\end{theorem*}

\begin{theorem*}[Corollary \ref{cor:manycomplements}]
There exist infinitely many subsets $\calS_0,\calS_1,\ldots$ of the set $\calP$ of primes,
all of lower density $0$,
all computable uniformly from an $\HTP(\Q)$-oracle (so that the rings $R_j=\Z[\calS^{-1}_j]$
are also uniformly computable below $\HTP(\Q)$),
with $\cup_j \calS_j=\calP$ and $\calS_i\cap \calS_j=\emptyset$ for all $i<j$,
and such that $\HTP(\Z_{\calS_j})\equiv_T \HTP(\Q)$ for every $j$.
\end{theorem*}

All of the above results generalize to {\bf all} number fields, which is
very different from the theorems that deal with diophantine models of $\Z$.
There the unconditional results apply only to number fields with extra
properties, such as the existence of an elliptic curve defined over $\Q$ and of rank one over $\Q$ and over the number field in question.

\section{Some facts from computability theory}
\label{S: definitions}
In this section we collected the definitions and facts from computability theory used in this paper.  We start with the most fundamental notions.
\begin{definition}A subset $T$ of $\Z$ is {\em computable}
(or {\em recursive}, or \emph{decidable}) if
there exists an algorithm (formally, a Turing machine;
informally, a computer program) that takes as
input any integer $t$ and, within finitely many steps, outputs YES or NO according to whether $t \in
T$.  The set $T$ is {\em computably enumerable} (or \emph{c.e.}, also
known as {\em recursively enumerable})
if it can be listed algorithmically:  some program,
running forever, outputs all the elements of $T$ (and nothing else),
although not necessarily in increasing order.
It is well known that there exist c.e.\ sets that are not computable.
\end{definition} 
Next we define the rings that are the focus of this paper.
\begin{definition}
A ring $R$ is \textit{computably presentable} if there is a
bijection $R \to \Z_{>0}$ such that the addition
and multiplication in $R$ correspond to computable functions 
$\Z_{>0}\times\Z_{>0} \to \Z_{>0}$ under this bijection.
The isomorphic copy of $R$ with domain $\Z_{>0}$ is said to be
a \emph{computable presentation} of $R$.
\end{definition}

 The next definition discusses the notion of an oracle used to define the Turing reducibility above.
\begin{definition}
Given a set $B \subseteq \Z$, we sometimes treat $B$ as an \emph{oracle}
for a computation of a function $f$.  This means that the computation follows a program,
but the program is allowed to use a subroutine which takes an input $n$
and returns the value of the characteristic function $\chi_B(n)$.  The oracle
set $B$ itself may not be computable; this is simply a method of saying 
that if we could compute $B$, then we would be able to compute $f$.
Such a function is said to be \emph{$B$-computable}.  A set $C$ is
$B$-computable (or \emph{Turing-reducible to $B$}, written $C\leq_T B$)
if $\chi_C$ is $B$-computable.
\end{definition}

For $\calS\subseteq\mathcal P$, computable presentability of the ring
$R=\Z[\calS^{-1}]$ is equivalent to computable enumerability of
$\calS$, and is also equivalent to computable enumerability of $\HTP(R)$,
where we consider $\HTP(R)$ as a subset of the positive integers in the
following sense.  Fix a computable enumeration of all polynomials (in
any number of variables) with integer coefficients and let $\HTP(R)$
be the set of indices corresponding to the polynomials with a root in
$R^k$ for the appropriate $k$.   Assuming our enumeration of polynomials is computable in the sense that we have Turing programs which, given an index $e$, compute the total degree of $f_e$ and all the coefficients of $f_e$, we can always recover the polynomial from its index and given a polynomial determine its index in a uniformly computable manner.  Thus, as a matter of convenience,  we can  view $\HTP(R)$ as a set of polynomials, writing ``$f_e\in\HTP(R)$'' in place of ``$e \in HTP(R)$''.    (An elementary but useful introduction
to computable rings and fields appears in \cite{M08}.)

There are countably many programs (or algorithms) and they can be
listed effectively, i.e., there is an algorithm which, given $n \in \Z_{>0}$, 
determines the corresponding program.  Fix such an enumeration.
The  {\em Halting Problem} is the
c.e.\ set of pairs $(m,n)$ such that the $m$-th program terminates on
input $n$.  A classical fact from computability theory is that
every c.e.\ set is Turing reducible to the Halting Problem and that the Halting
Problem is not computable. As we have mentioned already, there are
infinitely many Turing degrees of computably enumerable sets
that are undecidable, but not Turing equivalent to the Halting
Problem: that is, an oracle for the Halting Problem
can decide membership in these sets but not the other way around.


\newpage
\section{Examples of rings $R$ with $\HTP(R) \equiv_T\HTP(\Q)$ }
\subsection{Constructing rings with an $\HTP(\Q)$-oracle}
\label{sec:HTPQoracle}

In this section we first present a simpler version of the construction
given in \S\ref{sec:computable}.  While it illustrates the role of definability of
integrality at a prime in the construction, the ring $R$ produced by this
simple version is not necessarily computably presentable.  It
is computable only relative to the oracle for $\HTP(\Q)$, and here
$\HTP(R)$ is as defined in \S\ref{sec:intro}.

\begin{prop}
\label{prop:oracle}
There is a subring $R=\Z[\calS^{-1}] \subseteq \Q$ with $\HTP(R) \equiv_T \HTP(\Q)$, where
$\calS \leq_T  \HTP(\Q)$, and $\calS$ can be chosen 
to be co-infinite as a subset of the primes.  
\end{prop}
\begin{proof}
Fix a computable enumeration $\la f_e\ra_{e\in\N}$ of
$\Z[X_1,X_2,\ldots]$. Let $\calS_0=\emptyset$, and ${\calU_0}= \emptyset$. We proceed in stages; the $n$th stage will determine $\calS_n$ and $\calU_n$ (which will both be finite).

Assume we have just completed stage $n \geq 0$.  Now consider the polynomial
equation $f_{n}(\mathbf X)=0$ and use the $\HTP(\Q)$-oracle together with Corollary \ref{corollary:uniformsemilocal} to determine whether this
polynomial equation has solutions in $\Z[\calP-\calU_n]$.  This is possible by Proposition~\ref{proposition:uniformsemilocal} below.  If the answer is ``no,''
then $f_{n}$ is put on the list of polynomials without solutions in
our ring.  If the answer is ``yes,'' then we add $f_{n}$ to the
list of polynomials with solutions in our ring and search for a
solution integral at primes in ${\calU_{n}}$.  Once we locate the
solution, we add all the primes which appear in the denominators of
this solution to $\calS_n$ to form $\calS_{n+1}$.  Finally, we set
$\calU_{n+1}=\calU_n\cup\{ p\}$, where $p$ is
the least prime not in $\calU_n\cup\calS_{n+1}$.
This completes stage $n+1$, with $\calU_{n+1}\cap\calS_{n+1}=\emptyset$.
We let $\calS=\bigcup_{n=0}^{\infty}\calS_n$ and
$\calU=\bigcup_{n=0}^{\infty}\calU_n$.

To determine whether the $n$-th prime $p_n$ is inverted in our ring, we
just need to follow the construction until at most the $(n+1)$-st stage,
since (by induction on $n$) $p_n$ must lie in $\calS_{n+1}\cup\calU_{n+1}$.
(Here we regard $2$ as $p_0$.)  At the same time, to determine whether
a given polynomial has solutions in our ring, all we
need to do is again wait for the stage of the construction where this
polynomial was processed.  Thus, we have $\calS \leq_T \HTP(\Q)$ and
HTP($\Z[\calS^{-1}]) \leq_T \HTP(\Q)$.  Finally, $\calU$ must be
infinite, since each $\calU_n$ contains $n$ primes, and $\calU=\overline{\calS}$,
since the $n$-th prime lies in the (disjoint) union $\calS_{n+1}\cup\calU_{n+1}$.
Thus $\calS$ is co-infinite. 
\end{proof}

\begin{remark}
In the following proof, it is possible to ensure that $\calS$ is infinite as well as co-infinite,
by adding a prime to $\calS_{n+1}$ each time a prime is added to $\calU_{n+1}$.
If the process in the proof happened to build a finite $\calS$,
then we have a stronger statement: in this case, $\HTP(\Q)$ would have to compute
$\HTP(\Z)$, since $\HTP(\Z)\equiv_T\HTP(\Z[\calS^{-1}])$ for all finite $\calS$.
\end{remark}

Even though $\HTP(\Q)$ is c.e.\ and the set $\calS$ in 
Proposition \ref{prop:oracle} is enumerated using only an $\HTP(\Q)$-oracle,
$\calS$ itself need not be c.e.  The process of enumerating $\calS$ required
asking membership questions about $\HTP(\Q)$, and the computable
enumeration of $\HTP(\Q)$ is not sufficient to answer such questions.
Therefore, the subring $R$ built here need not be computably presentable.
(Of course, if $\HTP(\Q)$ should turn out to be decidable, then the procedure in
the proposition would indeed enumerate $\calS$ computably.)

\subsection{Constructing computably presentable rings with a priority argument}
\label{sec:computable}

From now on, we let $\calP$ denote the set of all rational primes,
let $\calW= \{q_1, q_2, q_3, \ldots\} \subseteq \calP$ be an infinite c.e.\ set
of primes, written in order of enumeration (not necessarily in increasing order),
and let $\calW_s = \{q_1, \ldots, q_s\}$.  Further, let
$\calR \subseteq \calP - \calW$ be any computable set.  Let $\calM=\calW\cup \calR$.  Observe that $\calM$ is still c.e.  The reader is
encouraged to assume that $\calR = \emptyset$ (and therefore $\calM=\calW$) for the first reading.  In
Theorem \ref{thm:low} and Theorem \ref{thm:B} we will use a nonempty $\calR$.

We construct a computably presentable ring $R$ of the form
$\Z[\calS^{-1}]$ satisfying $\HTP(R) \equiv_T \HTP(\Z[\calM^{-1}])$,
with $\calW \cup \calR=\calM$ and with $\calS$ being a computably
enumerable infinite and co-infinite set in $\calM$.  (In Theorem
\ref{thm:low}, we will have $\calR\subset \calS$, but in Theorem
\ref{thm:B} only a proper subset of $\calR$ will be in $\calS$ to be inverted.)
In particular, if $\calW = \calP$ and $\calR =\emptyset$, then we get
a computably presentable ring $R$ such that $\HTP(R) \equiv_T \HTP(\Q)$.
We remind the reader that the preceding section did not accomplish this, since the ring $R$ constructed
there required the use of an $\HTP(\Q)$-oracle in its construction.
In contrast, the ring $R$ that we build here is entirely defined
by an effective algorithm, with no oracle required; this makes it a
computably presentable ring as defined in the introduction.

 The rings of the form
$\Z[\calS^{-1}]$ with $\calS \subseteq \calP$ computably enumerable (but not
necessarily co-infinite) are precisely the computably presentable
subrings of $\Q$ (which are not necessarily computable but necessarily c.e).  We let
\[
{\tt E}:= \{(f, \vec x, j) : \exists n\geq 0\text{~such that~} f \in \Z[X_1, \ldots, X_n], \vec x \in
(\Z[\calM^{-1}])^n, j \in \Z_{>0}\},
\]
\label{ttE}
and we let $g\colon \Z_{>0} \to \tt E$ be a computable bijection.  We
do not actually make explicit use of the last coordinate $j$; its role
is to ensure that the pair $(f, \vec x)$ appears in the sequence
infinitely often.

\subsubsection{The construction of the set $\calS$}
\label{SS: constructionS}

Let $\{f_e \}_{e \in \Z_{>0}}$ be an enumeration of all polynomials $f
\in \Z[X_1, \ldots]$.  For each $e\geq 0$, we introduce a boolean
variable ${\tt R}_e$, to be updated depending on whether $f_e
\in \HTP(\Z[\calS^{-1}])$ in the course of the construction.   At the
beginning of the construction, ${\tt R}_e$ is set to FALSE for all $e
> 0$, and ${\tt R}_0$ is set to TRUE at the beginning of the construction. In our construction, at each step $s$, we will define sets
$\calS_s$ and $\calV_s$; at the end of the construction, we will
define $\calS = \cup_s \calS_s$. Let $h: \Z_{\geq0} \to \Z_{\geq 0}$ be a
computable function satisfying $h(0) = 0$  and several other conditions to be specified later.

\begin{itemize}
\item At stage $s = 0$: let $\calS_0 = \emptyset, \calV_0 = \emptyset$, {\tt R}$_0=$TRUE.
\item At stage $s > 0$: Let $e, \vec x,j$ be such that $g(s) = (f_e, \vec x, j)$.
Here $f_e$ is the $e$-th polynomial in the enumeration
$\{f_1, f_2, \ldots\}$ defined above.  Let
$\calV_{s}$ be the set of the first $h(s)$ elements in the enumeration of $(\calW-\calS_{s-1})$,  obtained from the given enumeration of $\calW$ by removing elements of $\calS_{s-1}$. 
Since $\calR\cap\calW=\emptyset$, we also have $\calR\cap\calV_s=\emptyset$.

\begin{itemize}
\item[(1)] If $\tt R_e =\textup{FALSE}$, $\vec x \in
  \Z[(\calM - \calV_s)^{-1}]^n$ and $f_e(\vec x) = 0$, then let
  $\calS_s := \calS_{s-1} \cup \calT_s$, where $\calT_s$ 
is the set of primes used in the denominators in $\xvec$.
Notice that $\calS_s\cap\calV_s=\emptyset$,
since none of these primes lies in $\calV_s$.
Here, we set $\tt R_e = \textup{TRUE}$.
\item[(2)] Otherwise, let $\calS_{s} := \calS_{s-1}$.
\end{itemize}
\item In the end, let $\calS := \cup_{s=1}^{\infty}\calS_s$, and $R = \Z[\calS^{-1}]$.  Since this construction
is entirely effective, $\calS$ is computably enumerable, and so $R$ is computably presentable.
\end{itemize}
\begin{remark}
\label{R: construction}
$\left.\right.$
\begin{itemize}
\item[(1)]
Every element of $\calR$ must belong to $\calS$, as $\calV_s$ never
stops elements of $\calR$ from being added to $\calS$. For each
$p \in \calR$, at the stage $s$ with $g(s) = (px-1, \frac 1 p, 1)$, we will be in Case (1)
of the above construction, and so $p$ will be added to the set $\calS$ of primes to be
inverted if $p$ was not already inverted. 
\item[(2)] 
Checking whether  $\vec x \in \Z[(\calM - \calV_s)^{-1}]^n$ is an
algorithmic operation since we are given that  $\vec x \in
  \Z[\calM^{-1}]^n$.  Thus we just need to check that none of the finitely many primes of $\calV_s$ occur in the denominator of $\vec{x}$.

\end{itemize}
\end{remark}

\begin{definition}
Let $s_0 = 0$. For each $e \geq 1$, we define $s_e$ to be the smallest positive integer strictly bigger than $s_{e-1}$ such that all of
the eventually true variables ${\tt R}_1, \ldots, {\tt R}_{e-1}$ have
been set to \textup{TRUE} by stage $s_e$ of the above construction,
and such that $g(s_e)=(f_e, \vec x, j)$.
\end{definition}

The stage $s_e$ must exist, since only finitely many
  ${\tt R}_i$ are considered to define it.  Computing $s_e$ requires
  an oracle, as we will see in Theorem~\ref{thm:originalpriority}.
In the following sections, we will show that $\HTP(\Z[\calM^{-1}])
\geq_T \HTP(R)$, given that the function $h$ satisfies the
following conditions:
\begin{align}
\label{E: condition1} 
&\textup{$h(s_e) \to
  \infty$ as $e \to \infty$}.\\
\label{E: condition2}
&\textup{$h(s_e) \leq h(s')$ when $s' > s_{e}$ satisfies $g(s') = (f_{e'}, \vec x, j)$ for some $e'$, $\vec x$, and $j$,
with $e' > e$.}\\
%
\label{E: condition3} 
&h(s) = h(s_e)\textup{ whenever }g(s) = (f_e,\vec x,j) \textup{ for some }e,\vec x, j\textup{ with } s\geq s_e.
\end{align}
\begin{remark}
Our definition of the stages $s_e$ and condition (2) imply that $h(s_{e-1}) \leq h(s_e)$, since $s_{e-1} <s_e$ and $g(s_{e-1})=(f_{e-1},\vec{x},j)$ and $g(s_{e})=(f_{e},\vec{x}',j')$ for an appropriate $\vec x, \vec x', j, j'$.
\end{remark}
For example, the function $h$ given by $h(s) = e$ (where $g(s)=(f_e, \vec x, j))$
clearly satisfies this set of conditions (since $h(s)=h(s_e)=e$ for
any stage $s$ at which $g(s)=(f_e, \vec x, j)$). It may make it
easier to understand the paper
to assume this for now, and also that  $\calW = \calP$. 
 In Corollary~\ref{C: easyh} we will use the function
$h(s)=e$. However, later on in \S\ref{generalh}, we will make other
choices of $h$ and $\calW$ as well.

\subsubsection{Analyzing the complement of $\calS$ in $\calW$}
\label{SS: complement}
 We now describe a
procedure that determines the complement of $\calS$ in order of enumeration in $\calW$.  In
Theorem \ref{thm:originalpriority} we show that the complement can be
computed with the aid of an oracle for $\HTP(\Z[\calM^{-1}])$.  As $\calS$ satisfies $\calR \subseteq \calS \subseteq \calR \cup \calW=\calM$, we focus on describing the set of primes in $\calW$ that are not in $\calS$, which we denote by $\calW - \calS$.
\begin{prop}\label{complement}
For $e>0$, the first $h(s_e)$ primes in $\calW-\calS_{s_e}$ 
(under our enumeration of $\calW$) are precisely
the first $h(s_e)$ primes in $\calW-\calS$.   
\end{prop}
%
%
%
%

\begin{proof}
  

It suffices to prove for any fixed $e>0$ that, at every stage $t \geq s_e$, the first $h(s_e)$
primes in $\calW -\calS_{s_e}$, which we denote by $p_1,p_2,\ldots,p_{h(s_e)}$
in order of enumeration into $\calW$, are the first $h(s_e)$ primes in $\calW - \calS_t$. 

We use induction on $t$.  For the base case $t = s_e$, the claim is
true by definition of the set $\calW - \calS_{s_e}$. For the induction
step, assume $\{p_1,\ldots, p_{h(s_e)}\}$ are the first $h(s_e)$
elements of $\calW - \calS_{t-1}$ (in the order of enumeration of
$\calW$).  Now consider the stage $t$ of Construction \ref{SS:
  constructionS} with $g(t) = (f_{e'}, \vec x, j).$ If we are in Case
(2) at this stage, then $\calS_{t-1} = \calS_{t}$ and
$\{p_1,\ldots, p_{h(s_e)}\}$ are the first $h(s_e)$ elements of
$\calW - \calS_{t}$.

 If we are in Case (1), then we must have $e' \geq e$, since otherwise
 ${\tt R}_{e'}$ would change from FALSE to TRUE at stage $t>s_e$,
 contradicting the definition of $s_e$.  Also, by Conditions (2) and
 (3) we have that $h(t) \geq h(s_e)$, with Construction \ref{SS:
   constructionS} implying that $\calV_t$ (the set of primes that
 cannot be inverted at this stage) is the set of the first $h(t)$
 elements of $\calW \setminus \calS_{t-1}$.  Thus,
 $\{p_1,\ldots, p_{h(s_e)}\} \subset \calV_t $ while
 $\calV_t\cap\calT_t=\emptyset$. We remind the reader that $\calT_t$
 is the set of primes in the denominator of a root of $f_{e'}$ being
 processed at this stage $t$. Also, by definition,
 $\calS_t =\calS_{t-1}\cup\calT_t$.  Thus,
 $\{p_1, \ldots, p_{h(s_e)}\} \subseteq \calW- \calS_t$.
\end{proof}

\begin{remark}
Since $h(s_e) \to
  \infty$ as $e \to \infty$,
it follows that the set $\calW -\calS$ is infinite.

\end{remark}
We now show that an $\HTP(\Z[\calM^{-1}])$-oracle is enough to determine polynomials with
solutions in the constructed ring $R = \Z[\calS^{-1}]$.  The proof is divided into two
parts.  First, we show that the membership of $f_e$ in $\HTP(R)$ can be determined by the end of stage $s_e$ (Proposition~\ref{P: equivalence}). Then we show that the $s_e$ can be computed using an oracle for $\HTP(\Z[\calM^{-1}])$.
\begin{prop}
\label{P: equivalence}
 For each $e >0$, the following are equivalent:
\begin{itemize}
\item[(i)] ${\tt R}_e = \textup{TRUE}$ at some stage $s$ (thus, at all stages $\geq s$);
\item[(ii)] $f_e \in \HTP(R)$;
\item[(iii)]  $f_e \in \HTP(\Z[(\calM - {\{p_1, \ldots, p_{h(s_e)}\}})^{-1}])$.
\end{itemize}
\end{prop}
\begin{proof}
The implication (i)$\Rightarrow$(ii) is trivial, the implication (ii)$\Rightarrow$(iii) follows directly from Proposition \ref{complement}, and so we need only to prove (iii)$\Rightarrow$(i). Suppose
  that $f_e \in \HTP(\Z[(\calM-{\{p_1, \ldots, p_{h(s_e)}\}})^{-1}])$. This
  means that there exists a solution $\vec x \in \Z[(\calM-{\{p_1, \ldots, p_{h(s_e)}\}})^{-1}]$
  such that $f_e(\vec x) = 0$. Consider a stage
  $s >s_{e}$ in  Construction \ref{SS: constructionS} such that $g(s) = (f_e, \vec x, j)$ for some $j$.
   
   If ${\tt R}_e$ is already set to TRUE by this stage, we are done. Otherwise, we must be in Case (1) at this stage of the construction, since
 \begin{align*}
  \calV _s &= \{\textup{first $h(s)$ primes of $\calW$ not in $\calS_{s-1}$}\} \\
  &= \{\textup{first $h(s_e)$ primes of $\calW$ not in $\calS_{s-1}$} \} \quad \textup{(Condition (\ref{E: condition3}))}\\
&= \{p_1, \ldots, p_{h(s_e)}\} \quad \textup{(Proposition~\ref{complement})}.
  \end{align*}
Then $\vec x \in \HTP(\Z[(\calM-{\{p_1, \ldots, p_{h(s_e)}\}})^{-1}]) = \HTP(\Z[(\calM-\calV_s)^{-1}])$, so
  ${\tt R}_e$ is set to \textup{TRUE} at this stage.
  \end{proof}
  
  This proposition serves as an important tool for computing of $s_e$ using the $\HTP(\Z[\calM^{-1}])$-oracle via Corollary \ref{corollary:uniformsemilocal} and finitely many primes in place
  of the whole set $\calS$.

\begin{theorem}
\label{thm:originalpriority}
The ring $R = \Z[\calS^{-1}]$ satisfies $\HTP(\Z[\calM^{-1}]) \geq_T \HTP(R)$.
\end{theorem}
\begin{proof}
  To show that $\HTP(\Z[\calM^{-1}]) \geq_T \HTP(R)$ it is enough to
  show that given $e \in \Z_{>0}$, the oracle for
  $\HTP(\Z[\calM^{-1}])$ can compute $s_e$.  Indeed, we show that
  knowing $s_e$ is enough to determine whether $f_e \in \HTP(R)$.  By
  Proposition \ref{P: equivalence}, we have that $f_e \in \HTP(R)$ if
  and only if
  $f_e \in \HTP(\Z[(\calM - {\{p_1, \ldots, p_{h(s_e)}\}})^{-1}])$.
  If we can compute $s_e$, we can run Construction \ref{SS:
    constructionS} for $s_e$ stages so that we can determine the
  primes $p_1$, \ldots, $p_{h(s_e)}$.  Now these primes can be plugged
  into the procedure discussed in Corollary
  \ref{corollary:uniformsemilocal} to determine whether
  $f_e \in \HTP(\Z[(\calM - {\{p_1, \ldots, p_{h(s_e)}\}})^{-1}])$
  using the oracle for $\HTP(\Z[\calM^{-1}])$.

  We now show how, given an $e \in \Z_{>0}$, the oracle for $\HTP(\Z[\calM^{-1}])$ computes $s_e$.  We proceed by induction.  \\
  {\it Base case}.  By definition, $s_1$ is the smallest $s>0$ such
  that $g(s)=(f_1,\vec x, j)$ for some $\vec x$ and some $j$.  This
  $s$ is certainly computable by running Construction \ref{SS:
    constructionS} through a sufficient number of stages until the
  first component of $g(s)$ is equal to $f_1$. Also, such a stage
  exists because
  {\it every} triple of the form $(f_1, \vec x, j)$ is in the
  image of $g(s)$.)

\vspace{.4em}
\noindent{\it Induction step.}  Suppose that using the oracle for
$\HTP(\Z[\calM^{-1}])$ we have computed $s_e$.  In this case, without
loss of generality, we can assume that we have also computed
$p_1, \ldots, p_{h(s_e)}$, and therefore can determine whether
$f_e \in \HTP(R)$, as described above.  If $f_e \not \in \HTP(R)$,
then $s_{e+1}$ is the smallest $s >s_e$ such that
$g(s)=(f_{e+1},\vec x,j)$ for some $\vec x$ and some $j$.  Such an $s$
can certainly computed by running Construction \ref{SS: constructionS}
through a sufficient number of stages.  If $f_e \in \HTP(R)$, then we
first run Construction \ref{SS: constructionS} through a sufficient
number of stages to find the smallest possible stage $\bar s \geq s_e$
such that $g(\bar s)=(e,\vec x, j)$ for some $\vec x$ and some $j$
with $f_e(\vec x)=0$.  The existence of $\bar s$ is guaranteed by the
fact that $f_e \in \HTP(R)$ and each pair $(f_e, \vec x)$ appears with
infinitely many $j$.  Now we look for the smallest $s >\bar s$ such
that $g(s)=(f_{e+1}, \vec x, j)$ for some $\vec x$ and some $j$.  We
set $s_{e+1}=s$.
\end{proof}

\begin{corollary}
\label{C: easyh}
There exists a subset $\calS \subseteq \calP$ that is infinite and co-infinite, such that the ring $R[\calS^{-1}]$ satisfies $\HTP(\Q) \equiv_T \HTP(R)$.
\end{corollary}
\begin{proof}
Choose $h(s)= e$, where $g(s) = (f_e, \vec x, j)$, $\calR = \emptyset$, and $\calW=\calP$.
\end{proof}


\subsection{Imposing Density Conditions}\label{generalh}
In this section, we show that a more complicated choice of $h$ lets us impose a density condition on the set $\calS$. We first need to define the relative upper and lower density of sets of primes. 
Given a subset $\calS
\subseteq \calW=\{q_{1},\ldots, q_{n}, \ldots\}$, we define
\[
\limsup_{n\rightarrow \infty} \frac{|\calS \cap \{q_{1},\ldots, q_{n}\}|}{n},
\]
to be the {\em upper density of a subset $\calS$ relative to $\calW$}.
The {\em relative lower density} can be defined in a similar
fashion.  Both of these depend on the choice of the enumeration of $\calW$,
as well as on $\calW$ itself:  a set $\calS$ could have upper density
$0$ in $\calW$ under one enumeration of $\calW$, yet have upper
density $1$ in the same set $\calW$ under a different enumeration.
However, if $\calW = \calP$ and we enumerate the primes
into $\calP$ in increasing order, then the upper
(resp. lower) density of $\calS$ relative to $\calW$ matches the usual
notion of upper (resp. lower) density. If the upper and lower
density are equal, this quantity is the {\em natural density} of $\calS$ in $\calW$.

We will use the following notation for our density
calculations in the following sections.
\begin{notation}
\label{notation:density}
 Suppose that $\calU \subset \calW$ is finite. Let $i_{\calW}(\calU)$ be the largest index $i$ such that $q_i \in \calU$
(or $1$ if $\calU$ is empty). We define
\[
 \mu_{\calW}(\calU):=\frac{|\calU|}{i_{\calW}(\calU) }.
 \]
\end{notation}
\begin{remark}[Constructing sets with $\mu_{\calW}$ arbitrarily close to 1, containing a given finite set $\calU$ and avoiding a given finite set $\calS$]
\label{rem:mu}
Let $\calS \subset \calW$ be a finite set  such that $\calS\cap \calU=\emptyset$.  Let $i=i_{\calW}(\calU \cup \calS)$ and consider a set $\calU_k=\calU\cup\{q_{i+1},q_{i+2}, \ldots, q_{i+k}\}$.  Now $\mu_{\calW}(\calU_k)=\frac{|\calU_k|}{i_{\calW}(\calU_k) }=\frac{ |\calU | + k}{i+k} \longrightarrow 1$ as $k \longrightarrow \infty$ and at the same time $\calU_k \cap \calS =\emptyset$.
\end{remark}
\begin{remark}
\label{first}
Let $\calS \subset \calW$ and let $t \in \Z_{>0}$.  Let $\calV$ be the
set of the first $t$ elements of $\calW$ not in $\calS$, and let
$\calV'$ be a set of any $h$ elements not in $\calS$.  We have that
$i_{\calW}(\calV) \leq i_{\calW}(\calV')$and therefore
$\mu_{\calW}(\calV) \geq \mu_{\calW}(\calV')$.
\end{remark}

Define a function $h(s)$ to be equal to the least positive integer greater or equal to $e$ that yields
$\mu_{\calW}(\calV_s) > \frac{e}{e+1}$ (where $e$ is such that $g(s) = (f_e, \vec x, j)$).
With this choice $h(s)$ is computable and well-defined by Remark \ref{rem:mu} because $\calS_{s-1}$ is a finite set.

\begin{prop}
\label{P: condition}
The function $h(s)$ satisfies Conditions~(\ref{E: condition1})-(\ref{E: condition3}).
\end{prop}
\begin{proof}
First, Condition~(\ref{E: condition1}) is trivially satisfied, since by definition of $s_e$ we have that $g(s_e)=(f_e,\vec x, j)$ for some $\vec x$ and some $j$, and therefore $h(s_e) \geq e$ and $h(s_e)\rightarrow \infty$ as $e \rightarrow \infty$.

To show that Condition~ (\ref{E: condition2}) is satisfied, we first
show that if $e\geq e'$ and $s \geq s'$ while $g(s)=(f_e, \vec x, j),
g(s')=(f_{e'},\vec x', j')$, then $h(s) \geq h(s')$.  First of all,
$\calS_{s'-1} \subseteq \calS_{s-1}$ and $\frac{e}{e+1} \geq
\frac{e'}{e'+1}$.  Now assume $h(s) <h(s')$.  Now we have, by the definition of $h(s)$,
\begin{equation}
\label{minimality}
e' \leq e \leq h(s) <h(s').
\end{equation}
By Remark \ref{first}, if $\calV'_{s}$ were any set of the form $\{\textup{\textit{any} $h(s)$ primes not contained in $\calS_{s-1}$}\}$, then $\mu_{\calW}(\calV_{s}) \geq \mu_{\calW}(\calV_{s}')$.
The primes not contained in $\calS_{s-1}$ are also not contained in $\calS_{s'-1}$
(although the $h(s)$ primes above are not necessarily the first $h(s)$ ones in the enumeration of $\calW$), and so
\begin{align*}
\mu_{\calW}(\calV_{s}) &:= \mu_{\calW}(\textup{first $h(s)$ primes not in $\calS_{s-1}$})\\
& = \mu_{\calW}(\textup{the same $h(s)$ primes as above not in $\calS_{s'-1}$})\\
& \leq \mu_{\calW}(\textup{first $h(s)$ primes not in $\calS_{s'-1}$})(\mbox{by Remark \ref{first}})\\
& \leq \frac{e'}{e'+1}\textup{ (by minimality of the choice of $h(s')$ and \eqref{minimality})} \\
& < \mu_{\calW}(\textup{the first $h(s')$ primes not in $\calS_{s'-1}$}) =: \mu_{\calW}(\calV_{s'}).
\end{align*}
But then $\frac{e'}{e'+1} \geq \mu_{\calW}(\calV_{s}) > \frac{e}{e+1}$, which is a contradiction when $e \geq e'$.  Thus Condition~ (\ref{E: condition2}) is satisfied.

To show that  Condition~ (\ref{E: condition3}) holds we need to use an induction argument similar to the one used in Proposition \ref{complement}.  (We cannot use the Proposition \ref{complement} itself as it was proved under the assumption that $h(s)$ satisfied the three conditions under consideration.) The argument is almost identical.  We again need to consider what happens to the primes of $\calV_{s_e}$, which are the first $h(s_e)$ non-inverted primes of $\calW$ at the end of the stage $s_e$, between the stage $s_e$ and any other stage $\tilde s$ with $g(\tilde s)=(f_e, \vec x, j)$ for some $\vec x$ and some $j$.  More specifically, we want to show that these primes do not get inverted.  As in Proposition \ref{complement}, we need to consider two cases: we are at the stage $s'>s_e$ with the corresponding $e' < e$ or at the stage $s'>s_e$ with the corresponding $e' \geq e$.  In the first case, by definition of $s_e$, we must be in the second clause of Construction \ref{SS: constructionS} and no new primes are inverted.  In the second case, $h(s') \geq h(s_e)$ by Condition~ (\ref{E: condition2}), and so the only primes which can be inverted at this stage are primes not in $\calV_{s_e}$.  

Thus, when the construction starts the stage $\tilde s$, the primes of $\calV_{s_e}$ are not inverted and $\calV_{s_e}$ must satisfy the requirements for $h(\tilde s)$ that are the same as the requirements for $h(s_e)$.  So by the minimality requirement on $h(s)$, we must have $h(s_e)=h(\tilde s)$.
\end{proof}

\begin{theorem}
\label{T: density}
For each infinite c.e.\ set $\calW$ of primes and each computable enumeration of $\calW$,
there exists a set $\calS \subseteq \calW$ such that $\calW - \calS$ has relative upper density $1$ in $\calW$, and such that $\HTP(\Z[\calW^{-1}]) \geq_T \HTP(\Z[\calS^{-1}])$. In particular, $\calS$ has relative lower density $0$ in $\calW$.
\end{theorem}
\begin{proof}
Set $\calR = \emptyset$ in the above construction, so that using notation of  Theorem~\ref{thm:originalpriority} we have $\calW=\calM$. That $\HTP(\Z[\calW^{-1}]) \geq_T \HTP(\Z[\calS^{-1}])$ now follows from Theorem~\ref{thm:originalpriority}. As for the density, we consider the stages $s_e$ for $e \in \Z_{>0}$, at which we have
\[
\mu_{\calW}(\{p_{i_1},\ldots,p_{i_{h(s_e)}}\}) > \frac{e}{e+1}.
\]
Considering these stages $s_e$, we see that $\lim_{e \rightarrow
  \infty}\mu_{\calW}(\calV_{s_e})=1$, so that the relative upper density of $\calS$ in $\calW$ is $1$.
\end{proof}

From the proposition above we can obtain the following corollary, which requires repeated applications of the proposition.  
\begin{corollary}
\label{cor:nested}
There exists a sequence $\calP=\calW_0 \supset \calW_1\supset \calW_2
\ldots$ of uniformly c.e.\ sets of rational primes (with $\calP$ denoting the
set of all primes, enumerated in ascending order) such that
\begin{enumerate}
\item  $\HTP(\Z[\calW_i^{-1}]) \equiv_T \HTP(\Q)$ for $i \in \Z_{>0}$, 
\item $\calW_{i-1} - \calW_i$ has the relative upper density (with respect to
the enumeration of $\calW_{i-1}$) equal to 1 for all $i \in \Z_{>0}$,
\item The relative lower density of $\calW_i$ (with respect to the enumeration
of $\calW_{i-1}$) is 0, for $i \in \Z_{>0}$.
\end{enumerate}
\end{corollary}
\begin{proof}
By transitivity we can arrange that $\HTP(\Q) \geq_T
\HTP(\Z[\calW_i^{-1}])$ for all $i$.  On the other hand,
Corollary~\ref{cor:RQ} implies that $\HTP(\Q) \leq_T \HTP(\Z[\calW_i^{-1}])$.  
\end{proof}

\subsection{Arbitrary lower density}
As remarked in the introduction, it is desirable to have
the density for the set of inverted primes to be equal to 0.
However, it seems that this is difficult to
accomplish using this construction.   We can, however, control the lower density:
\begin{theorem}
\label{thm:low}
For every computable real number $r$ between 0 and 1 there exists a
c.e.\ set $\calS$ of primes such that the lower density of $\calS$ is $r$
and $\HTP(\Z[\calS^{-1}]) \equiv_T \HTP(\Q)$.
 \end{theorem}
To be clear:  here, a \emph{computable real number} $r$ is a real
number such that some computable sequence of rational numbers
has limit $r$.  The upper and lower cuts defined in $\Q$ by $r$
will therefore be Turing-reducible to the Halting Problem $\emptyset'$,
but need not be decidable.  This definition, common in number
theory, is less strict than the usual meaning of the same phrase
in computability theory.
\begin{proof}
  Let $r$ be a computable real number. We set the parameters of the
  construction of \S\ref{SS: constructionS} as follows: $\calM=\calP$;
  $\calR \subseteq \calP$ is some computable set of density $r$ (to understand why such a set of primes always exists, see 
  \cite{EES});
  $\calW := \calP-\calR$, enumerated in ascending order (computably,
  since $\calR$ is computable);
  and $h(s)$ the least integer $\geq e$ such that
  $\mu_{\calW}(\calV_s) > \frac{e}{e+1}$, where $g(s)=(f_e,\xvec,j)$,
   exactly as in the construction in \ref{SS: constructionS}.  (This $h(s)$ also
  satisfies Conditions~(\ref{E: condition1})-(\ref{E: condition3}) for
  the same reason as Proposition~\ref{P: condition}).

Then as in Theorem~\ref{T: density}, the set $\calS_{\calW}:=\calS \cap (\calP - \calR)=\calS\cap\calW$ has relative lower density 0 in $\calW$.

Further, this choice of parameters ensures that all the primes of
$\calR$ are inverted by the end of the construction, as in
Remark~\ref{R: construction}.  It remains to show that the set
$\calR \cup \calS_{\calW}=\calS$ is of (absolute) lower density $r$.
Given that the set $\calS_{\calW}$ is of lower relative density $0$,
the set $\calW$ is enumerated in increasing order and the set $\calR$
is of absolute density $r$, there exists a positive integer sequence
$\{m_i\} \rightarrow \infty$ with the property that for every
$\varepsilon >0$ there exists $M>0$ such that for all $i>M$
\begin{equation}
\label{density:1}
 \frac{\#\{q \in \calS_{\calW}, q \leq m_i\}}{\#\{q \in \calW, q \leq m_i\}} +\frac{\#\{q \in \calR, q \leq m_i\}}{\#\{q \in \calP, q \leq m_i\}} \leq \varepsilon.
 \end{equation}
 At the same time,
 \begin{equation}
 \label{density:2}
 \frac{\#\{q \in \calS, q \leq m_i\}}{\#\{q \in \calP, q \leq m_i\}}=\frac{\#\{q \in \calS_{\calW}, q \leq m_i\}}{\#\{q \in \calP, q \leq m_i\}} +\frac{\#\{q \in \calR, q \leq m_i\}}{\#\{q \in \calP, q \leq m_i\}}
 \end{equation}
and for all $i$,
\begin{equation}
\label{density:3}
\frac{\#\{q \in \calS_{\calW}, q \leq m_i\}}{\#\{q \in \calP, q \leq m_i\}} \leq \frac{\#\{q \in \calS_{\calW}, q \leq m_i\}}{\#\{q \in \calW, q \leq m_i\}}.
\end{equation}
Thus, the desired density assertion follows from \eqref{density:1} -- \eqref{density:3}.
\end{proof}

\subsection{Complementary Rings}
\label{sec:density}

We also show that by modifying the construction 
of Section \ref{sec:HTPQoracle} to obtain several rings $R_i=\Z[\calS_i^{-1}]$ at once,
we may arrange that the sets $\calS_i$ form a partition of the set of all primes. 
In this case, the sets $\calS_i$ are not built to be computably enumerable,
and so the rings $R_i$ may not be computably presentable.
\begin{theorem}
\label{thm:complements}
For any positive integer $m$, the set of all rational primes $\calP$
can be represented as a union of pairwise disjoint sets $\calS_1,
\ldots, \calS_{m}$, each of upper density 1 and such that for all $i$
we have that $\HTP(\Z[\calS_i^{-1}] )\leq_T \HTP(\Q)$ and $\calS_i
\leq_T \HTP(\Q)$.
\end{theorem}
\begin{proof}
  We imitate the construction given in Section \ref{sec:HTPQoracle},
  defining all $\calS_{1,s},\ldots,\calS_{m,s}$ so that their pairwise
  intersections are always empty, but their union is an initial
  segment of $\Z$, growing ever larger as $n$ increases.  (As before
  we will have $\calS_i=\bigcup_{s=1}^{\infty}\calS_{i,s}$.)  Start
  with $\calS_{i,0}=\emptyset$ for all $i=1,\ldots,m$ and fix the same
  computable list $\la f_e\ra_{e\in\Z_{>0}}$ as before.  Also, set all
  ${\tt R}_{e,i}, i=1,\ldots,m$ to FALSE.  We proceed in stages with
  $i=0,\ldots, m-1$.  For stage $s=me+i >0$, let
  $\calW_s=\bigcup_{j\not =i}\calS_{j,s-1}$, and use the oracle for
  $\HTP(\Q)$ to determine whether $f_e$ has a root $\vec x$ in
  $\Z[\overline \calW_s^{-1}]$.  We go through a different sequence of
  actions depending on the answer.
\begin{itemize}
\item[``Yes'':] Let $\calT_{i,s}$ be the smallest set of primes such
  that $\vec x \in \Z[\calT_{i,s}^{-1}]$.  Next set
  ${\tt R}_{i,e}=TRUE$, and let $i_s$ be the index (in the listing of
  all primes) of the largest prime used in the construction so far.
  Now set
 \[
\calS_{i,s}=\calS_{i,s-1} \cup \calT_{i,s} \cup \{ \mbox{every prime with index in the interval }[1,2^{i_s}] \mbox{ not in } \calW_s~\},
\]
and set $\calS_{j,s}=\calS_{j,s-1}$ for $j\not =i$.  (The reason for
including the last set of primes is to ensure upper density 1 for
$\calS_{i}$.  Indeed, in the computation of the upper density, 
at this point we get from
$\calS_{i,s}$ a term that will be greater than or equal to
$\frac{2^{i_s}}{i_s+2^{i_s}}$, so it converges to $1$ as
$s \rightarrow \infty$.)
\item[``No'':] As above, let $i_s$ be the index (in the listing of all primes) of the largest prime used in the construction so far.  Now set
 \[
\calS_{i,s}=\calS_{i,s-1} \cup \{ \mbox{every prime with index in the interval }[1,2^{i_s}] \mbox{ not in } \calW_s ~\},
\]
and set $\calS_{j,s}=\calS_{j,s-1}$ for $j\not =i$.
\end{itemize}
To determine whether $f_e$ has a solution in $\Z[\calS_i^{-1}]$, all
we need to do is run the construction until the step $me+i$.  If ${\tt
  R}_{e,i}$ is not set to TRUE st this point, it never will be.
Similarly, to determine if a prime $p_j$ is in $\calS_i$ we just need
for the construction to process this prime to see where it was put.
Thus, both assertions concerning Turing reducibility are true.
\end{proof}

One can also build countably many rings with the same property, instead of building $m$ rings as in Theorem~\ref{thm:complements},

\begin{corollary}
\label{cor:manycomplements}
There exist infinitely many subsets $\calS_1,\calS_2,\ldots$ of the set $\calP$ of primes,
all of lower density $0$,
all computable uniformly from an $\HTP(\Q)$-oracle (so that the rings $R_j=\Z[\calS^{-1}_j]$
are also uniformly computable below $\HTP(\Q)$),
with the property that $\cup_j \calS_j=\calP$ and $\calS_i\cap \calS_j=\emptyset$ for all $i<j$,
and such that $\HTP(\Z[\calS_j^{-1}])\equiv_T \HTP(\Q)$ for every $j$.
\end{corollary}
\begin{proof}
  Use the construction from Theorem \ref{thm:complements}, but loop
  through a computable listing of all pairs $(i,e)\in\N^2$ with the first index
  referring to the set of primes and the second to a polynomial.
\end{proof}
\begin{remark}
\label{rmk:computablesubring}
Notice that, if the sets $\calS_i$ in Theorem \ref{thm:complements}
were all c.e.\ (or if those in Corollary \ref{cor:manycomplements}
were uniformly c.e.), then they would all be computable.
It is not clear that the theorems of this section fail to hold when one demands
that the set(s) of inverted primes be computable, but the constructions
given here do not suffice. (In Remark \ref{rmk:low}, we will
mention how close we can come.) It is an open question, worthy of attention,
whether there exists a computable coinfinite set $\calS$ of primes
with $\HTP(\Z[\calS^{-1}])\equiv_T\HTP(\Q)$.
\end{remark}

\section{Rings $R$ such that $R \equiv_T \HTP(R) \equiv_T \HTP(\Q)$} 
\label{sec:generalize}
In this section we construct examples where $\HTP$ of a ring is no more complicated than
the ring itself.  One such example is already known:  if $\calS$ is a computably
enumerable subset of $\calP$ with Turing degree $\bfz'$, such as the Halting Problem itself,
then $\HTP(\Z[\calS^{-1}])$ is also computably enumerable, hence $\leq_T\calS$.
However, we wish to build sets $\calS$ of Turing degree $<\bfz'$:  we will show
that every computably enumerable Turing degree $\geq_T\HTP(\Q)$ contains
a c.e.\ set $\calS$ with $\HTP(\Z[\calS^{-1}])\equiv_T\calS$.  (Of course,
if $\HTP(\Q)$ is decidable, then $\Q$ would be another example of this phenomenon.)
\begin{theorem}
\label{thm:B}
For every computably enumerable set $B \subset \Z_{>0}$ with
$\HTP(\Q)\leq_T B$, there exists a computably presentable ring
$\Z[\calS^{-1}]$, with $\calS$ a computably enumerable subset of the
prime numbers of lower density 0, such that $\calS\equiv_T
\HTP(\Z[\calS^{-1}])\equiv_T B$.  (By setting $B\equiv_T \HTP(\Q)$, we obtain a ring
$R$ with $\HTP(R)\equiv_T \HTP(\Q)$.)
\end{theorem}
We prove the theorem in the remainder of this section. The required construction, while similar to Construction \ref{SS: constructionS}, is more complicated.  The proof proceeds as follows: first we construct a computably
enumerable set $\calS$. Then we prove that $\calS$ possesses certain properties
that are used to show that $B \geq_T \HTP(\Z[\calS^{-1}])$ (Theorem~\ref{thm:high}) and $\calS\geq_T B$ (Proposition~\ref{prop:high}).  This of course will imply that $\calS \equiv_T  \HTP(\Z[\calS^{-1}])$, since $\HTP(\Z[\calS^{-1}]) \geq_T \calS$ for obvious reasons.

Before we get started  we review and introduce some notation.
\begin{itemize}
\item  $H: \Z_{>0} \longrightarrow B$ will denote a computable function enumerating $B$.
\item $\calR$ is a computable infinite sequence of primes of density  0.  As before, let $\calW=\calP-\calR$ and observe that $\calW$ is computable.
\item For each $s>1$, let $y_s \in \calR$ be  the $x$-th least element
  of $\calR -\calS_{s}$ where $x=H(s)$.  
\end{itemize}
To ensure that $B$ can be computed from $\calS$ we will arrange that the following is true:

{\bf $(*)$}      $\forall x, s, t \in \Z_{>0}$,  if   $H(s) =x $ and the sets $\calR- \calS$ and $\calR - \calS_t$ agree on their least  $x$ elements, then  $t > s$.  

\smallskip

Given a positive integer $x$,  let 
\[
t_x=\min\{t\in \Z_{>0}|\calR- \calS \mbox{ and }\calR - \calS_t \mbox{ agree on their least  }x \mbox{ elements}\}.
\]
  It is not hard to see that $t_x$ is computable
from the $\calS$-oracle (as a function of $x$).  Indeed,  all we have to do is to run the construction and
compare the first $x$ elements of $\calR- \calS$ and
$\calR - \calS_i$, where $i$ is the current stage.  Eventually the
first $x$ elements of the two sets under consideration will be the
same.  The earliest stage at which this happens is $t_x$.  Now if $(*)$ holds we know that $t_x$ is
greater than $s$, the stage (if one exists) such that $H(s)=x$.

Therefore, assuming $(*)$ holds, to determine if $x\in B$, we just
have to check values of $H$ from $H(1)$ up to $H(t_x)$ to see if any
of them is equal to $x$.

To make sure that $t_x >s$, the $x$-th least element of
$\calR - \calS_s$ (denoted above by $y_s$) will be added to $\calS_{s+1}$ so that $\calR - \calS_s$ and $\calR - \calS$ cannot agree on their least $x$ elements.  

To see why the set differences cannot agree on the least $x$ elements any earlier, note the following, more general phenomenon:
if an integer $z \in (\calR - \calS_t) \cap (\calR- \calS)$, then $z \in \calR- \calS_u$ for all $u \geq t$, since we only add things to sets $\calS_u$, never remove them. Therefore, if for some $u <s$ we had that $\calR-\calS_u$ and $\calR-\calS$ agree on their least $x$ elements, say $a_1, \ldots, a_x$, then $a_1,\ldots,a_x \in\calR -\calS_s \subseteq \calR-\calS_u$ and $a_1,\ldots,a_x$ are still the least $x$ elements of the sets $\calR -\calS_{s}$ and $\calR -\calS$, contradicting that $a_x=y_s \in \calS_{s+1} \subseteq \calS$ in our construction.

\subsection{The construction of the set $\calS$}
\label{constructionhigh}
As we mentioned above, this construction is similar to Construction \ref{SS: constructionS},  but modified to
accommodate the equivalence $(*)$.  As before, for each $e \geq 1$, we
introduce a boolean variable ${\tt R}_e$ and we also have ${\tt R}_0$=TRUE at the beginning of the construction. At the beginning of the
construction, ${\tt R}_e$ is set to FALSE for all $e > 0$. At each
step $s$ of the construction, we define the sets $\calS_s$ and
$\calV_s$; at the end of the construction, we will define $\calS =
\cup_s \calS_s$.  Below the ``complement'' notation will denote the
complement in the set of all primes again.

\begin{itemize}
\item At stage $s = 0$: let $\calS_0 = \emptyset, \calV_0 = \emptyset$.
\item At stage $s > 0$: suppose $g(s) = (f_e, \vec x, j) \in \tt E$, 
where $f_e$ is the $e$-th polynomial in the enumeration
$\{f_1, f_2, \ldots\}$ as defined above.   At this stage,
we will consider this polynomial $f_e$.  

We now reconsider the function $h(s)$ and the set $\calV_s$.  The function $h(s)$ is  defined the same way as above in Proposition \ref{P: condition},  i.e., $h(s)$ is the smallest integer greater or equal to $e$ such that $\mu_{\calW}(\mbox{the set of the first $h(s)$ primes of $\calW$ not in } \calS_{s-1}) > \frac{e}{e+1}$.  

  Let 
\[
\calV_s =\{p_{i_1}, \ldots, p_{i_{h(s)}}\} \cup  \{ \mbox{ the least } e \mbox{ elements of }  (\calR- (\calS_{s-1} \cup \{y_s\})\}.
\]
\end{itemize}
Now there are two cases:

\begin{itemize}
\item[(1)] If ${\tt R}_e =\textup{FALSE}$,
  $\vec x \in \Z[(\calP - {\calV_s})^{-1}]$ and $f_e(\vec x) = 0$,
  then let $\calS_s = \calS_{s-1} \cup \calT_s \cup\{y_s\}$, where
  $\calT_s$ is the least set of primes such that
  $\vec x \in (\Z[\calT_s^{-1}])^n$.  (The set $\calT_s$ is allowed to
  have elements from $\calR - \calV_s$.) Here, we set
  ${\tt R}_e = \textup{TRUE}$.
\item[(2)] Otherwise, $\calS_s=\calS_{s-1} \cup\{y_s\}$.
\end{itemize}
In the end, let $\calS = \cup_{s=1}^{\infty}\calS_s$, and $R = \Z[\calS^{-1}]$.  Since this construction
is entirely effective, $\calS$ is computably enumerable, and so $R$ is computably presentable.

\subsection{Analyzing the complement of $\calS$} We now describe the
procedure determining the complement of $\calS$ in order.  The
following notation will be useful: Let $\calV := (\calP- \calS) \cap
(\calP- \calR)=\{p_{i_1}, p_{i_2}, \ldots\}$ in order of a fixed
enumeration, and let $\calH:=(\calP- \calS) \cap
\calR=\{w_1,w_2,\ldots\}$ in order of a fixed enumeration. We will use
the familiar variable $s_e$, but it will be defined inductively in a
slightly different manner, together with another variable $v_e$.  We
define
\begin{itemize}
\item $s_0=0$
\item For $e \geq 1$ we let $v_e$ be the smallest positive integer greater than $s_{e-1}$ such that
for each positive integer $i \in [1,e]$, either $i \not \in B$ or for some $s < v_e$ we
have that $H(s)=i$.
\item For $e \geq 1$ we let $s_e$ be the smallest possible positive integer greater than $v_e$ such that $g(s_e)=(f_e,\vec x, j)$ for some $\vec x $ and some $j$ and ${\tt R_{e-1}}$ has either been set to \textup{TRUE} or never will be.
\end{itemize}

We start with the following proposition.


\begin{prop}
At every stage $u \geq s_e$, we have that $\{p_{i_1} \ldots, p_{i_{h(s_e)}}\} \subset \calV_u$. 
\end{prop}
\begin{proof}
  This proof is analogous to the proofs of Proposition \ref{P: condition} and Proposition \ref{complement}, since the function $h(s)$ and the sequence $\{s_e\}$ satisfy the assumptions of these propositions.
  \end{proof}
To determine $\calH=\{w_1<w_2<\ldots\}$ in order we use the following proposition:
\begin{prop}
  If $\{w_1,\ldots, w_e\} = \mbox{ the least } e \mbox{ elements of
  } \calR \cap  \calV_{s_e}$, then for every
  $t \geq 0$, we have that $\{w_1,\ldots, w_e\} \subset
  \calV_{t+s_e}$.
\end{prop}
\begin{proof}
  We use strong induction on $e$ again.  The base case of $t=0$ is
  clear.  So assume that $t >0$ and $\{w_1,\ldots, w_e\} = \mbox{
    the least } e \mbox{ elements of }\calR\cap
  \calV_{s_e-1+t}$.  Let $g(t+s_e)=(f_{e'}, \vec x, j) \in {\tt
    E}$ (as defined on page \pageref{ttE}).  If $e' \leq e$, then $\calT_{t+s_e}=\emptyset$ (since ${\tt
    R}_{e'} $ has either been satisfied already or will never be).
  Thus, only $y_{t+s_e}$ is added to $\calS_{t-1+s_e}$ to form
  $\calS_{t+s_e}$.  By assumption on $v_e < s_e$, we have that $H(t+s_e) >
  e$, and therefore the first $e$ elements of $\calR -
  \calS_{t-1+s_e}$ and $\calR - \calS_{t+s_e}$ are the same.

Suppose now that $e'>e$.  In this case $\calV_{t+s_e}$ will contain $w_1,\ldots, w_e$ and $\calT_{t+s_e}$ cannot contain any of $w_1,\ldots,w_e$.  Further, $y_{t+s_e}$ cannot be among $w_1,\ldots,w_e$ for the same reason as above.
\end{proof}
\begin{prop}
\label{P: equivalencehigh}
For each $e \geq 0$, the following are equivalent:
\begin{itemize}
\item[(i)] ${\tt R}_e = \textup{TRUE}$;
\item[(ii)] $f_e \in \HTP(\Z[\calS^{-1}])$;
\item[(iii)] $f_e \in \HTP(\Z[\left (\calP-{\calV_{s_e}}\right)^{-1}])$.
\end{itemize}
\end{prop}
\begin{proof}
The proof is analogous to the proof of Proposition~\ref{P: equivalence}.
\end{proof}

\begin{theorem}
\label{thm:high}
 $B \geq_T \HTP(\Z[\calS^{-1}])$.
\end{theorem}
\begin{proof}
We follow Construction \ref{constructionhigh}, using a $B$-oracle to compute stages $s_e$ (and hence the primes of $\calV_{s_e}$) recursively, starting with $s_0=0$.  Assume inductively that
we have computed $s_{e-1}$.  First, using the $B$-oracle, we check whether $e \in B$.
If $e\notin B$, then we know $v_e=s_{e-1}+1$;
while if $e\in B$, then either $v_e=1+s_{e-1}$ (if some $s\leq s_{e-1}$ has $H(s)=e$);
or else $v_e$ is the unique integer $v$ with $H(v)=e$ (if this $v$ is $>s_{e-1}$).
Since $H$ is computable, we have thus determined $v_e$.

Knowing $v_e$, we then find the smallest integer $s>v_e$ such that $g(s)=(f_{e-1},\vec x, j) \in \tt E$.  Compute $\calV_s$ for this stage and use the $B$-oracle to determine whether $f_{e-1} \in \HTP(\Z[(\calP-\calV_s)^{-1}])$.  (This is possible by Proposition \ref{proposition:uniformsemilocal},
because $\HTP(\Q)\leq_T B$.)
If the answer is ``no,'' then set $s_e=s$.  If the answer is ``yes,''
find the least stage $\hat s\geq s$ such that $g(\hat s)=(f_{e-1},\vec x, j')$, $\vec x \in \Z[(\calP-{ \calV_{\hat s}})^{-1}]=\Z[(\calP-{\calV_{s}})^{-1}]$, and $f_{e-1}(\vec x)=0$.  Then $s_e$ will be the least integer $\geq \hat s$ with
$g(s_e)=(f_e,\xvec,j)$ for some $\xvec$ and some $j$.
\end{proof}

The last proposition we need to complete the proof of Theorem \ref{thm:B} states the following.
\begin{prop}\label{prop:high}
$\calS \geq_T B$.
\end{prop}
\begin{proof}
We use an $\calS$-oracle to determine membership in $B$.   Given a
positive integer $x$, run the construction until it reaches a stage $s$ at which $\calP-
\calS$ and $\calP- \calS_s$ agree on at least their smallest $x$
elements.  If $x\notin\{ H(1),H(2),\ldots,H(s)\}$, then $x \not \in
B$ by property $(*)$, which was introduced right before Section~\ref{constructionhigh}.
Otherwise $x\in\rg{H}=B$.
\end{proof}
\begin{remark}[Density]
Our final observation concerns the assertion on density in Theorem \ref{thm:B}.  The construction above did not consider density, but the assertion can be fulfilled by applying the same methods as in Proposition \ref{T: density}.
\end{remark}

\begin{remark}
\label{rmk:low}
Mathematicians familiar with priority arguments will not be surprised to learn that standard
lowness requirements (e.g.\ from \cite[VII.1]{S87}) can be mixed in with certain of the constructions
in the preceding two sections. In Theorem \ref{thm:low}, for example,
we can build a c.e.\ set $\calS\subseteq\calP$ which is low (i.e., $\calS'\leq_T\emptyset'$),
has arbitrary computable lower density $r$, and satisfies $\HTP(\Q)\equiv_T\HTP(\Z[\calS^{-1}])$.
It follows that the entire subring
$\Z[\calS^{-1}]$ has low Turing degree as a subset of $\Q$. This is the closest
we currently come to answering the question in Remark \ref{rmk:computablesubring}.
The construction is a technical extension of that for Theorem \ref{thm:low}, and we do not include it here.

It is open whether a low set $\calS$ could satisfy Theorem \ref{thm:B},
where a c.e.\ set $B\geq_T\HTP(\Q)$ is given arbitrarily.
The requirement that $\calS\equiv_T B$ would have to be dropped, of course,
but it seems possible that a low c.e.\ set $\calS$, of lower density $0$,
might satisfy $\HTP(\Z[\calS^{-1}])\equiv_T B$. However, the techniques
we have used do not suffice to prove this.
\end{remark}

\section{Computing $\HTP(R)$ from $\HTP(\Q)$ for a semilocal ring $R \subseteq \Q$}
\label{sec:semilocal}

We begin with results necessary to show that $\HTP(R) \leq_T \HTP(\Q)$ for all rings $R \subseteq \Q$.

\begin{prop}[\cite{Sh5}, {Theorem~4.2}]
\label{P: nonzero}
  Let $R$ be any subring of $\Q$.  Then the set of non-zero elements
  of $R$ is existentially definable over $R$.
\end{prop}
This proposition  implies:
\begin{corollary}
\label{cor:RQ}
For every ring $R \subset \Q$, we have that $\HTP(\Q) \leq_T \HTP(R)$. 
\end{corollary}
\begin{proof}
  Every element of $\Q$ can be written
  as a ratio $x/y$ with $x,y \in R, y \not=0$.  Thus we can replace all the
  variables ranging over $\Q$ by quotients of variables ranging over
  $R$ and add equations stipulating that the denominators are not
  zero using Proposition~\ref{P: nonzero}.
\end{proof}

In a similar fashion one can show that the following corollary is true
as well.
\begin{corollary}
\label{relativized}
For every ring $R \subset \Q$, every $m>0$, and every set $A \subset \Q^m$
diophantine over $\Q$, the set $A \cap R^m$ is diophantine over $R$.
\end{corollary}

This corollary will be used in this section to prove that
an oracle for $\HTP(\Z[\calM^{-1}])$ also gives an oracle for $\HTP(\Z[(\calM-\calP_0)^{-1}])$,
for any finite
set of primes $\calP_0 \subseteq \calM$ (Corollary~\ref{corollary:uniformsemilocal}).

Now, we give an effective version of Julia Robinson's
result concerning the diophantine definition of elements in $\Q$ that
are integral at a finite set of primes. This lets us decide from an
$\HTP(\Q)$-oracle whether polynomials have solutions in certain semilocal rings
$\Z[(\calP-{\{p_1, \ldots, p_n\}})^{-1}]$.  We show that the following result is true:

\begin{prop}
\label{proposition:uniformsemilocal}
Let $\calP_0 = \{p_1, \ldots, p_n\}$ be any finite set of primes. Then
$\HTP(\Z[(\calP-\calP_0)^{-1}])$ is computable uniformly in
$\HTP(\Q)$ in the following sense: there exists an algorithm, using $\HTP(\Q)$ as an oracle,
that can decide, given any such $\calP_0$, whether any given polynomial with coefficients
in $\Q$ has a solution in $\Z[(\calP-\calP_0)^{-1}]$.
\end{prop}

\begin{lemma}
\label{L: }
Let
\[
f_{a,b}(x_{11}, x_{12}, x_{13}, x_{14}) = x_{11}^2-ax_{12}^2 - bx_{13}^2 + abx_{14}^2-1
\]
be a polynomial. Then for $i \in \{1,2\}$ and $j \in \{1, 2,3,4\}$:
\begin{itemize}
\item[(1)] 
\begin{align*}
\Z_{(2)} = \{2(x_{11} + x_{12} + y_{11} + y_{12}) | \exists x_{ij}, y_{ij} \in \Q: 
f_{3,3}(x_{1j})^2 + f_{3,3}(x_{2j})^2 + f_{2,5}(y_{1j})^2 + f_{2,5}(y_{2j})^2 = 0 \}.
\end{align*}
\item[(2)] If $p \equiv 3 \bmod 8$, then
\[
\Z_{(p)} = \{2(x_{11} + x_{12} + y_{11} + y_{12}) | \exists x_{ij}, y_{ij} \in \Q:
 f_{-1,p}(x_{1j})^2 + f_{-1,p}(x_{2j})^2 + f_{2,p}(y_{1j})^2 + f_{2,p}(y_{2j})^2 = 0 \}.
 \]
\item[(3)] If $p \equiv 5 \bmod 8$, then
\begin{align*}
\Z_{(p)} = \{&2(x_{11} + x_{12} + y_{11} + y_{12}) | \exists x_{ij}, y_{ij} \in \Q:\\
&f_{-2p,-p}(x_{1j})^2 + f_{-2p,-p}(x_{2j})^2 + f_{2p,-p}(y_{1j})^2 + f_{2p,-p}(y_{2j})^2 = 0 \}.
\end{align*}
\item[(4)] If $p \equiv 7 \bmod 8$, then
\[
\Z_{(p)} = \{2(x_{11} + x_{12} + y_{11} + y_{12}) | \exists x_{ij}, y_{ij} \in \Q:
 f_{-1,-p}(x_{1j})^2 + f_{-1,-p}(x_{2j})^2 + f_{2p,p}(y_{1j})^2 + f_{2p,p}(y_{2j})^2 = 0 \}.
\]
\item[(5)] If $p \equiv 1 \bmod 8$, and if $q$ is a prime congruent to $3$ mod $8$ such that $\left(\frac{p}{q}\right) = -1$, then
\[
\Z_{(p)} = \{2(x_{11} + x_{12} + y_{11} + y_{12}) | \exists x_{ij}, y_{ij} \in \Q:
 f_{-2p,q}(x_{1j})^2 + f_{-2p,q}(x_{2j})^2 + f_{2p,q}(y_{1j})^2 + f_{2p,q}(y_{2j})^2 = 0 \}.
\]
\end{itemize}
\end{lemma}
\begin{proof}
We first prove part (a) in detail. The set
\[
S_{a,b}:=\{2x_1 \in \Q | \exists x_2, x_3, x_4 \in \Q: f_{a,b}(x_1, x_2, x_3, x_4) = 0\}
\]
is equal to the set of norm-$1$ elements of the twisted quaternion
algebra
$\Q \cdot 1 \oplus \Q \cdot \alpha \oplus \Q \cdot \beta \oplus \Q
\cdot \alpha \beta$
with $\alpha^2 = a, \beta^2 = b$, and $\alpha\beta = -\beta\alpha$. By
\cite[Proposition~10(a)]{Koenig2}, we know that
\[
\Z_{(2)} = (S_{3,3} + S_{3,3}) + (S_{2,5} + S_{2,5}),
\]
from which the conclusion follows.

The rest of the lemma follows by similar methods of proof: parts (2)
to (4) are \cite[Proposition~10(b)]{Koenig2}, and part (5) is
\cite[Proposition~10(c)]{Koenig2}.  \qed\end{proof}

For simplicity in notation, let $q \in \{1,2,3,5,7\}$ satisfy
$q \equiv p \bmod 8$, and let $f_q$ be the polynomial that defines
$\Z_{(p)}$ as in the previous lemma. Now we can prove Proposition
\ref{proposition:uniformsemilocal}:

\begin{proof}[Proof of Proposition \ref{proposition:uniformsemilocal}]
  Let $\calP_0 = \{p_1, \ldots, p_n\}$, and let
  $Q(Z_1,\ldots,Z_k) \in \Z[Z_1,\ldots,Z_k]$. We need to ascertain
  whether $Q(Z_1,\ldots,Z_K) =0$ has any solutions in
  $\Z[\overline{\calP_0}^{-1}]$.  We proceed in the following fashion:
\begin{enumerate}
\item If $p_i$ is congruent to $1$ modulo $8$, choose a prime $q_i$
  such that $q_i \equiv 3 \bmod 8$ and
  $\left(\frac{p_i}{q_i}\right) = -1$.
\item Consider the following system, which can be regarded as a single polynomial using sums of squares:
\begin{equation}
\label{system}
\left \{
\begin{array}{c}
Q(Z_1,\ldots,Z_k) =0\\
f_{(p_i \bmod 8)}(Z_i, x_{i2},x_{i3},x_{i4})=0, 1 \leq i \leq n\\
\end{array}
\right .
\end{equation}
\end{enumerate}
The above system has rational solutions $\{Z_i,x_{ij}\}$ if and only
if $Q(Z_1,\ldots,Z_k) =0$ has rational solutions integral at
$p_1,\ldots, p_n$, by the construction of the polynomials $f_{\ell}$,
$\ell \in \{1,2,3,5,7\}$.
\end{proof}

We also need Proposition~\ref{proposition:uniformsemilocal} in a more general setting, since the
general case of our
construction in \S\ref{sec:computable} gives rings $R$ such that  $\HTP(R) \equiv_T
\HTP(\Z[\calM^{-1}])$.
\begin{corollary}
\label{corollary:uniformsemilocal} Let $\calM$  be any set of
primes and let
$\calP_0 = \{p_1, \ldots, p_n\}\subseteq\calM$ be any finite set of primes
from $\calM$. Then
$\HTP(\Z[(\calM-\calP_0)^{-1}])$ is computable uniformly in
$\HTP(\Z[\calM^{-1}])$ in the following sense: there exists an algorithm,
using $\HTP(\Z[\calM^{-1}])$ as an oracle, that can decide, given any $\calP_0$,
whether any given polynomial with coefficients in $\Z[(\calM-\calP_0)^{-1}]$ has a solution in
$\Z[(\calM-\calP_0)^{-1}]$.
\end{corollary}

\begin{proof}
This is an immediate consequence of Corollary~\ref{relativized}.
\end{proof}

\section{Number Fields}
\label{sec:numberfields}
In this section, we discuss the extensions of our results to number
fields.  Throughout this section, $K$ will denote a number field, and
$\calO_K$ its ring of integers. Let $\calP_K$ denote the set of finite
primes of $K$, and if $\calS$ is a set of prime ideals, let $\calO_{K, \calS}$ denote the ring of
$\calS$-integers, which is defined by
\[
\calO_{K, \calS} = \{x \in K : \ord_{\pp}x \geq 0 \textup{ for all } \pp \not\in \calS\}.
\]

 We first give a brief survey of what is known about Hilbert's Tenth Problem over finite extensions of $\Q$.
\subsection{Hilbert's Tenth Problem for subrings of number fields}

Hilbert's Tenth Problem for the ring of integers of a number field has
been resolved in many cases, but undecidability for the ring of
integers of {\em any} number field is only known under the assumption that
the Shafarevich-Tate conjecture holds \cite{MR}. The theorem below summarizes
what is known.
\begin{theorem}
The ring  $\Z$ has a diophantine definition and Hilbert's Tenth Problem is
  undecidable over the rings of integers of the following fields:%
\begin{itemize}%
\item Extensions of degree 4, totally real number fields and their
  extensions of degree 2. (See \cite{Den2} and \cite{Den3}.)  These
  fields include all Abelian extensions.
\item Number fields with exactly one pair of non-real embeddings (See \cite{Ph1} and \cite{Sh2}.)%
\item Any number field $K$ such that there exists an elliptic curve $E$ of positive rank defined over $\Q$ with $[E(K):E(\Q)] < \infty$. (See \cite{Po} and \cite{Sh33}.)%
\item Any number field $K$ such that there exists an elliptic curve of rank 1 over $K$ and an Abelian variety over $\Q$
keeping its rank over $K$. (See \cite{CPZ}.)
\end{itemize}%
\end{theorem}%
\comment{
\begin{definition}
\label{def:bigrings}
For a number field $K$, let $\mathcal P_K$ denote the set of finite
primes of $K$, i.e.\ the set of all prime ideals of $O_{K}$, the ring of integers of $K$.  Given a
set~$\calS$ of prime ideals, not necessarily finite, the ring~$O_{K,\calS}$
is defined to be the subring of~$K$ defined by
$$O_{K,\calS}=\{x\in K:\ord_{\pp}x \geq 0  \mbox{ for all } \pp \notin \calS\}.$$
Observe that if $\calS= \emptyset$, then $O_{K,\calS}=O_K$ and if
$\calS=\calP_K$, then $O_{K,\calS}=K$.  If $\calS$ is finite, $O_{K,\calS}$ is
called a {\em ring of $\calS$-integers}.
\end{definition}
}
To measure the ``size'' of a set of primes of a number field one can
also use the natural density defined below for a number field.
\begin{definition}
\label{def:natural}
Let $\calS \subseteq \calP_K$.
The {\em natural density} of $\calS$ is defined to be the limit
\[
    \lim_{X\rightarrow \infty}
    \frac{\#\{\pp \in S: N\pp \leq X\}}
    {\#\{\pp \in \calP_K: N\pp \leq X\}}
\]
if it exists. (Here $N\pp$ denotes the size of the residue field of the prime or its norm.)  If the limit above does not exist, one can talk about
{\it upper} density by substituting $\limsup$ for $\lim$, or {\it
  lower} density by substituting $\liminf$ for $\lim$.
\end{definition}
 The proposition below summarizes what we know about $\calS$-integers of number fields.
 \begin{theorem}
$\left.\right.$
 \begin{itemize}
 \item If $K$ is a totally real number field, an
   extension of degree 2 of a totally real number field or such that
   there exists an elliptic curve defined over $\Q$ and of the same
   positive rank over $K$ and $\Q$, then for any $\varepsilon >0$,
   there exists a set ${\mathcal W}$ of primes of $K$ whose natural
   density is bigger than $1-[K:\mathbb{Q}]$$^{-1} - \varepsilon$ and
   such that $\mathbb{Z}$ has a diophantine definition over
   $O_{K,{\mathcal W}}$, thus implying that Hilbert's Tenth Problem is
   undecidable over $O_{K,{\mathcal W}}$. (See \cite{Sh1}, \cite{Sh3},
   \cite{Sh6} and \cite{Sh33}.)
\item Assume there is an elliptic curve defined over $K$ with $K$-rank
  equal to 1. For every $t>1$ and every collection
  $\delta_1,\ldots, \delta_t$ of nonnegative computable real numbers adding
  up to 1, the set of primes of $K$ may be
  partitioned into $t$ mutually disjoint computable subsets $\calS_1,\dots
  ,\calS_t$ of natural densities $\delta_1,\ldots, \delta_t,$ respectively,
  with the property that $\Z$ admits a diophantine model in each ring
  $O_{K,\calS_i}$. In particular, Hilbert's Tenth Problem is undecidable for
  each ring $O_{K,\calS_i}$. (See \cite{PS}, \cite{EE},\cite{Perlega}, \cite{EES}.)
 \end{itemize}
 \end{theorem}
In~\cite{MR}, Mazur and Rubin showed  that if the
 Shafarevich-Tate conjecture holds, then for every cyclic extension of
 number fields $M/K$
 of prime degree there always exists an elliptic curve defined over
 $K$
 with $K$-rank
 and $M$-rank
 equal to one.  So if the Shafarevich-Tate conjecture holds for all
 number fields, one can show the undecidability of Hilbert's
 Tenth Problem for the ring of integers and big rings for all number
 fields.

\subsection{Presenting primes of a number field in a computable manner}
Before proceeding with generalizations of our results, we need to
discuss how we are going to present primes of number fields, which,
unlike primes of $\Q$, do not necessarily correspond to a single number
but are ideals, i.e.\ infinite subsets of the field.  This data can be kept track of, since
rings of integers are Dedekind domains, which have finite presentations.

Given a prime $\pp$ of $K$, let $a(\pp) \in K$ be such that
  $\ord_{\pp}a(\pp)=-1$, $a(\pp)$ has non-negative order at all other
  primes of $K$, and for every $K$-prime $\qq \not = \pp$, conjugate to
  $\pp$ over $\Q$, we have that $\ord_{\qq}a(\pp)=0$.
The proof of the following proposition can be found in \cite{Bach}.
\begin{prop}
\label{computable}
Given a number field $K$, the following statements are true:
\begin{enumerate}
\item There exists a computable procedure which outputs the following
  information for each rational prime $p$: $(a(\pp_1),e_1, f_1,
  \ldots, a(\pp_m), e_m, f_m)$, where $p=\pp_1^{e_1}\ldots
  \pp_m^{e_m}$ is the factorization of $p$ (the prime ideal
  corresponding to $p$, really) in $K$ and $f_i$ is the relative
  degree of $\pp_i$ over $p$.
\item There exists a computable procedure that for each $x \in K$
  outputs the following information: $(q_1, a(\qq_1), \ord_{\qq_1}x,
  \ldots, q_m, a(\qq_m), \ord_{\qq_m}x)$, where $q_i$ is a rational
  prime number, $\qq_i$ is a prime of $K$, and $q_i =\qq_i \cap \Z$.
  Moreover, for each $K$-prime $\qq$ such that $a(\qq)$ is not on the list, we
  have that $\ord_{\qq}x=0$.
\item There exists a computable procedure that for each $x \in K$ determines whether all of the real conjugates of $x$ are positive.
\end{enumerate}
\end{prop}
Using the proposition above we identify a prime $\pp$ of a number
field $K$ with a pair $(p, a(\pp))$, where $p=\pp\cap \Z$.  We will fix a basis for $K$
over $\Q$ and represent each $a(\pp)$ by its rational coordinates with
respect to the chosen basis.  So finally a prime of $K$ will
correspond to an $(n+1)$-tuple of rational numbers, where $[K:\Q]=n$.
Given a set $K$-primes, we will now say that it is computable or
c.e.\ if the corresponding set of $n+1$-tuples of rational numbers is
computable or c.e.

\subsection{Effective diophantine definition of integrality at a finite set of primes.}

Let $K$ be a number field, and let $\pp_1, \ldots, \pp_n$ be a finite set
of primes of $K$. In this section, we would like to prove the analogue
of Proposition \ref{proposition:uniformsemilocal}. That is:

\begin{prop}
\label{proposition:uniformsemilocalnumberfield}
Let $\calP_0 = \{\pp_1, \ldots, \pp_n\}$ be any finite set of
primes. Let $R = \bigcap_{1 \leq i \leq n} \calO_{\pp_i} \cap K$. Then
$\HTP(R)$ is computable uniformly in $\HTP(K)$, in the sense of
Proposition \ref{proposition:uniformsemilocal}.
\end{prop}

Let $a, b \in K^{\times}$ be chosen as in \cite[Lemma~3.19]{Par13}. Then we have the Artin homomorphism that is defined on the prime ideals as 
\begin{align*}
\Psi: I^S &\to \textup{Gal}(K(\sqrt a, \sqrt b)/K) \cong \{ \pm 1 \}^2 \\
\pp &\mapsto \left(\left(\frac{a}{\pp}\right),\left(\frac{b}{\pp}\right)\right)
\end{align*}
and extended linearly to all ideals in $I^S$, where $S$ consists of all the infinite places and the places lying over $2ab$.

\begin{lemma}
\label{lemma:1,1}
There exists a polynomial $f_{(1,1)} \in K[t, y_1, y_2, x_1, \ldots, x_m]$ for some $m \geq 1$ such that for all prime ideals $\pp \in I^S$ satisfying $\Psi(\pp) = (1,1)$, there exists elements $y_{\pp,1}, y_{\pp,2} \in K$ such that
\[
\calO_{\pp} \cap K = \set{t \in K}{(\exists x_1, \ldots, x_m)~f_{(1,1)}(t, y_{\pp,1},y_{\pp,2}, x_1, \ldots, x_m) = 0}.
\]
\end{lemma}
\begin{proof}
This is \cite[Lemma 3.25(c)]{Par13}.
\end{proof}

\begin{lemma}
\label{lemma:not1,1}
Let $\sigma \in \{\pm 1\}^2$ with $\sigma \neq (1,1)$. Then there
exists a polynomial $f_{\sigma} \in K[t, y, x_1, \ldots,
x_{m_{\sigma}}]$ for some $m_{\sigma} \geq 1$ such that for all prime
ideals $\pp \in I^S$ satisfying $\Psi(\pp) = \sigma$, there exists an
element $y_{\pp} \in K$ such that
\[
\calO_{\pp} \cap K = \set{t \in K}{(\exists x_1, \ldots, x_m)~f_{(1,1)}(t, y_{\pp}, x_1, \ldots, x_m) = 0}.
\]
\end{lemma}
\begin{proof}
  With \cite[Proposition 2.3 and Definition 3.10]{Par13}, all that remains to
  do is to describe an algorithm that finds the element $y_{\pp}$ such
  that the Hilbert symbols satisfy $(y_{\pp}, a)_{\pp} =
  (y_{\pp},b)_{\pp} = -1$, but such that at least one of
  $(y_{\pp},a)_{\qq}$ or $(y_{\pp},b)_{\qq}$ is equal to $1$ for all
  places $\qq \neq \pp$. Such choice of $y$ exists due to \cite[ Theorem 3.7]{Par13}.
  Since there are explicit methods for computing Hilbert
  symbols, one enumerates all elements of $K$ and computes the Hilbert
  symbols $(y,a)_{\qq}$ and $(y,b)_{\qq}$ at the places $\qq$ that
  satisfy $v_{\qq}(y) \equiv 1 \bmod 2$, until we find such a $y$ (we
  only need to check finitely many places $\qq$ because of \cite[Lemma 3.8]{Par13}).
Enumerating through the $y$'s must terminate, because of
  the existence of such a $y$. The $f_{\sigma}$ can be described
  explicitly in terms of $y_{\pp}$, as in the previous subsection.
\end{proof}

\begin{lemma}
For any prime ideal $\pp$ of $\calO_K$, the localization ring $(\calO_K)_{\pp}$ is Diophantine.
\end{lemma}
\begin{proof}
This is \cite[Proposition 2.2]{PS}.
\end{proof}

Now the proof of Proposition
\ref{proposition:uniformsemilocalnumberfield} proceeds exactly like
that of Proposition \ref{proposition:uniformsemilocal}: Given a
polynomial $Q \in K[Z_1, \ldots, Z_k]$, we impose integrality
conditions on the $Z_{\ell}$ at each of the primes $\pp_i \in \calP_0$
using Lemma \ref{lemma:1,1} and Lemma \ref{lemma:not1,1}, as in
Proposition \ref{proposition:uniformsemilocal}.

\subsection{Turing reducibility between Hilbert's Tenth Problem for
  subrings of a number field $K$}
In this section we use Proposition
\ref{proposition:uniformsemilocalnumberfield} to prove results concerning
big subrings of number fields along the lines of similar results for
$\Q$.  We apply this proposition in the same manner as over $\Q$ to
construct rings $R$ with $\HTP(R) \equiv_T \HTP(K)$.

There are several well known and easy to prove Turing relations
between several types of rings and fields:
\begin{prop}
$\left.\right.$
\begin{enumerate}  
\item $\HTP(K) \leq_T \HTP(\Q)$
\item If $\calS$ is a cofinite set of $K$-primes, then $\HTP(\calO_{K,\calS}) \equiv_T \HTP(K)$.
\item For any set of $K$-primes $\calS$, including an empty set, we have $\HTP(K) \leq_T \HTP(\calO_{K,\calS})$.
\item For any set $\calS$ of $K$-primes  with a diophantine definition or model of $\Z$, we have $\HTP(\calO_{K,\calS}) \equiv_T \HTP(\Z)$. 
\item For any c.e.\ set of $K$-primes $\calS$, including an empty set,  we have $\HTP(\calO_{K,\calS}) \leq_T \HTP(\Z)$.
\end{enumerate}
\end{prop}
We now state our new results, for which we omit the proofs, since they
are identical to the proofs of the results for $\Q$. The upper and lower relative natural density can be defined in
an analogous manner for sets of primes of number fields, as it was for
sets of rational primes. It will, in general, also depend on the
underlying set as well as its ordering, as is the case for $\Q$.
\begin{theorem}
\label{thm:nestedK}
If $\calW$ is a c.e.\ set of primes of a number field $K$, then it
contains a c.e.\ subset $\calS$ such that the relative upper density
of $\calV=\calW - \calS$ is $1$ and
$\HTP(\calO_{K,\calW})\geq_T\HTP(\calO_{K,\calS})$.
\end{theorem}
\begin{corollary}
There exists a sequence $\calP=\calW_0 \supset \calW_1\supset \calW_2 \ldots$ of c.e.\ sets of  primes of a number field $K$ (with $\calP$ denoting the set of all primes of $K$) such that
\begin{enumerate}
\item  $\HTP(\calO_{K,\calW_i}) \equiv_T \HTP(K)\leq_T \HTP(\Q)$ for $i \in \Z_{>0}$, 
\item $\calW_{i-1} - \calW_i$ has the relative upper density (with respect to $\calW_{i-1}$) equal to 1 for all $i \in \Z_{>0}$,
\item The lower density of $\calW_i$ is 0, for all $i \in \Z_{>0}$.
\end{enumerate}
\end{corollary}
\begin{corollary}
There exists a computably enumerable subset $\calW$ of $K$-primes, of lower natural density $0$,
such that $\HTP(\Q) \geq_T \HTP(K) \equiv_T \HTP(\calO_{K,\calW})$. 
\end{corollary}

\begin{theorem}
For each computable  real number $r$ between 0 and 1 there is a c.e.\ set $\calS$ of primes of $K$ such that the lower density of $\calS$ is $r$ and $\HTP(\calO_{K,\calS}) \equiv_T \HTP(K)\leq_T\HTP(\Q)$.
 \end{theorem}

\begin{theorem}
For every computably enumerable set $B \subset \Z_{>0}$ with $\HTP(\Q)\leq_T B$,
there exists a computably presentable ring $R=\calO_{K,\calS}$, with $\calS$ a computably enumerable subset of  primes of $K$ of lower density 0, such that $\calS\equiv_T R \equiv_T \HTP(R)\equiv_T B$.
(By setting $B\equiv_T \HTP(\Q)$, we can thus make $\HTP(R)\equiv_T \HTP(\Q)$, of course.)
\end{theorem}
\begin{theorem}
  For every positive integer $m$, the set of all $K$- primes $\calP$ can
  be represented as a union of pairwise disjoint sets $\calS_1,
  \ldots, \calS_{m}$, each of upper density 1 and such that for all
  $i$ we have that $\HTP(\calO_{K,\calS_i} )\equiv_T \HTP(K) \leq_T
  \HTP(\Q)$ and $\calS_i \leq_T \HTP(K)$.
\end{theorem}
\begin{theorem}
There exist infinitely many subsets $\calS_1,\calS_2,\ldots$ of the set $\calP$ of $K$-primes,
all of lower density $0$,
all computable uniformly from an $\HTP(K)$-oracle (so that the rings $R_j=\calO_{K,\calS_j}$
are also uniformly computable below $\HTP(K)$),
which have $\cup_j \calS_j=\calP$ and $\calS_i\cap \calS_j=\emptyset$ for all $i<j$,
and such that $\HTP(R_j)\equiv_T \HTP(K)$ for every $j$.
\end{theorem}

\bibliographystyle{alpha}
\bibliography{mybib}%
\end{document}